\documentclass[12pt]{amsart}
\usepackage[all, cmtip]{xy}

\usepackage{times,amsfonts,amsmath,amstext,amsbsy,amssymb,
amsopn,amsthm,upref,eucal, mathrsfs}

\newcommand \la {\lambda}

\newcommand \Prob {{\mathbb P}}

\newcommand \Var {{\mathrm {Var}}}

\newcommand \Ai {{\mathrm {Ai}}}

\newcommand \Conf {{\mathrm {Conf}}}

\newcommand \zr {{z}^{\prime}}

\newcommand \zp {{\mathbb Z}^{\prime}}

\newtheorem{theorem}{Theorem}[section]

\newtheorem{lemma}[theorem]{Lemma}

\newtheorem{corollary}[theorem]{Corollary}
\newtheorem{proposition}[theorem]{Proposition}

\begin{document}
\title [Rigidity of Determinantal Point Processes]{Rigidity of Determinantal Point Processes with the Airy, the Bessel and the Gamma  Kernel}

\begin{abstract}

A point process is said to be {\it rigid} if for any bounded domain in the phase space, the number of particles in the domain is almost surely determined by the restriction of the configuration to the complement of our bounded domain. 
The main result of   this paper is that determinantal 
point processes with the Airy, the Bessel and the Gamma kernels are rigid. The proof follows the scheme of
Ghosh \cite{G-exp},  Ghosh and Peres \cite{GP}: the main step is the construction of a sequence of additive statistics with variance going to zero. 

\end{abstract}

\author{ Alexander I. Bufetov}
\address{Aix-Marseille Universit{\'e}, CNRS, Centrale Marseille, I2M, UMR 7373}
\address{ The Steklov Institute of Mathematics, Moscow}

\address{
The Institute for Information Transmission Problems, Moscow}\address{
National Research University Higher School of Economics, Moscow}

\maketitle

\section{Introduction.}	
	\subsection{Rigid Point Processes.}
	  Let $M$ be a complete separable metric space.  Recall that a {\it configuration} on $M$ is a purely atomic Radon measure on $M$; in other words,  a collection of particles considered without regard to order and not admitting accumulation points in $M$. The space $\text{Conf}(M)$ of configurations on $M$ is itself a complete separable metric space with respect to the vague topology on the space of Radon measures.  A point process on $M$ is by definition a Borel probability measure on $\text{Conf}(M)$. 
	
	Given a bounded subset $B \subset M$ and a configuration $X \in \text{Conf}(M),$ let $ \#_B(X)$ stand for the number of particles of $X$ lying in $B$. Given a Borel subset $C \subset M,$ we let $\mathcal{F}_C$ be the $\sigma$-algebra generated by all random variables of the form $\#_B, B\subset C.$ If $\mathbb{P}$ is a point process on $M,$ then we write $\mathcal{F}_C^{\mathbb{P}}$ for the $\mathbb{P}$-completion of $\mathcal{F}_C$.

The following definition of rigidity of a point process is due to Ghosh \cite{G-exp} (cf. also Ghosh and Peres \cite{GP}).

	{\bf Definition.} A point process $\mathbb{P}$ on M is called {\bf rigid} if for any bounded Borel subset $B \subset M$ the random variable $\#_B$ is  $\mathcal{F}_{M\backslash B}^{\mathbb{P}}$-measurable.

	Let $\mu$ be a $\sigma$-finite Borel probability measure on $\mathbb{R}$, and let $\Pi(x,y)$ be the kernel of a locally trace-class operator of orthogonal projection acting in $L_2(\mathbb{R}, \mu)$. 
Recall that the determinantal point process $\mathbb{P}_{\Pi}$ is a Borel probability measure on
$\Conf(\mathbb{R})$   defined by the condition that for any bounded measurable function $g$, for which $g-1$ is supported in a bounded set $B$,
we have
\begin{equation}
\label{eq1}
\mathbb{E}_{\mathbb{P}_{\Pi}}\Psi_g
=\det\biggl(1+(g-1)\Pi\chi_{B}\biggr).
\end{equation}
The Fredholm determinant in~\eqref{eq1} is well-defined since
$\Pi$ is locally of tracel class.
The equation (\ref{eq1}) determines the measure $\Prob_{\Pi}$ uniquely.
For any
pairwise disjoint bounded Borel sets $B_1,\dotsc,B_l\subset \mathbb{R}$
and any  $z_1,\dotsc,z_l\in {\mathbb C}$ from (\ref{eq1}) we  have
$$\mathbb{E}_{\mathbb{P}_{\Pi}}z_1^{\#_{B_1}}\dotsb z_l^{\#_{B_l}}
=\det\biggl(1+\sum\limits_{j=1}^l(z_j-1)\chi_{B_j}\Pi\chi_{\sqcup_i B_i}\biggr).$$
For further results and background on determinantal point processes, 
see e.g. \cite{BorOx},  \cite{HoughEtAl}, \cite{Lyons}, 
\cite{Lytvynov},  \cite{ShirTaka1}, \cite{ShirTaka2}, \cite{Soshnikov}.

We now formulate a sufficient condition for the rigidity of a determinantal point process on $\mathbb{R}$.
\begin{proposition}\label{gen-rig}
Let $U\subset {\mathbb R}$ be an open subset,  let $\mu$ be the Lebesgue measure on $U$, and let $\Pi(x,y)$ be a kernel yielding an operator of orthogonal projection 
acting in  $L_2({\mathbb R}, \mu)$.
Assume that there exists $\alpha\in (0, 1/2)$, $\varepsilon>0$, and, for any $R>0$, a constant $C(R)>0$ such that  the following holds:
\begin{enumerate}
\item if $|x|, |y|\geq R$, then 
$$
|\Pi(x,y)|\leq C(R)\cdot \frac{({x}/{y})^{\alpha}+ ({y}/{x})^{\alpha}}{|x-y|};
$$
\item if $|x|\leq R$, then  for all $y$ we have 
$$
\int\limits_{x:|x|\leq R}|\Pi(x,y)|^2d\mu(x)\leq \frac{C(R)}{1+y^{1+\varepsilon}}.
$$

\end{enumerate} 
Then the point process $\Prob_{\Pi}$ is rigid. 
\end{proposition}

As we shall see below, Proposition \ref{gen-rig} implies  rigidity for determinantal point 
processes with the Airy and the Bessel kernels; in the last subsection of the paper, we shall 
obtain a counterpart of Proposition \ref{gen-rig} for determinantal point processes with discrete phase space and, as its corollary, rigidity  for the determinantal point 
process with the Gamma kernel.

{\bf Remark.}
As far as I know, rigidity of point processes  first appears (under a different name) in the work of Holroyd and Soo \cite{hsoo}, who established, in particular, that the determinantal point process with the Bergman kernel 
is {\it not} rigid. For the sine-process, rigidity is due to Ghosh \cite{G-exp}. For the Ginibre ensemble, rigidity has been established by Ghosh and Peres \cite{GP}; 
see also Osada and Shirai \cite{OS}. 

\subsection{Additive Functionals and Rigidity.}

Given a bounded measurable function $f$ on $M$, we introduce the additive functional $S_f$ on $\mathrm{Conf}(M)$ by the formula 
$$
	S_f(X) = \sum\limits_{x \in X} f(x).
$$
We recall the sufficient condition for rigidity of a point process given by Ghosh \cite{G-exp}, Ghosh and Peres \cite{GP}.
\begin{proposition} [Ghosh \cite{G-exp}, Ghosh and Peres \cite{GP}]\label{g-p-r}
	Let $\mathbb{P}$ be a Borel probability measure on $\mathrm{Conf}(M)$.  Assume that for any $\varepsilon >0$ and any bounded subset $B \subset M$ there exists a bounded measurable function $f$ of bounded support such that
$f\equiv 1$ on $B$ and 
$\mathrm{Var}_{\mathbb{P}} S_f < \varepsilon$.
Then the measure $\mathbb{P}$ is rigid.
\end{proposition}
{\bf Proof.} For the reader's convenience, we recall the  elegant  short proof of Ghosh \cite{G-exp}, Ghosh and Peres \cite{GP}. 
 Let $B^{(n)}$ be an increasing sequence of nested bounded Borel sets exhausting $M$. Our assumptions and the Borel-Cantelli Lemma imply the existence of a sequence of bounded measurable function $f^{(n)}$ of bounded support, such that $f^{(n)}|_{B^{(n)}} \equiv 1$ and that for $\mathbb{P}$-almost every $X \in \mathrm{Conf}(M)$ we have 
$$
	 \lim\limits_{n \to \infty} S_{f^{(n)}}(X) - \mathbb{E}_{\mathbb{P}} S_{f^{(n)}} = 0.
$$
	Since, for any bounded $B$ and sufficiently large $n$, we have 
$$
	S_{f^{(n)}}(X) = \#_B (X) + S_{f^{(n)}\chi_{M\backslash B}}(X),
$$
we thus obtain the equality 
$$
	\#_B (X) = \lim\limits_{n \to \infty} (-S_{f^{(n)}\chi_{M\backslash B}}(X) + \mathbb{E}S_{f^{(n)}}),
$$
for $\mathbb{P}$-almost every $X$, and the rigidity of $\mathbb{P}$ is proved.

{\bf Remark.} In fact, to prove rigidity, it suffices that the function $f$ only satisfy the inequality $|f-1|<\varepsilon$ on $B$; the proof of the proposition becomes slightly more involved, but the result is still valid. 
\subsection{Variance of Additive Functionals.}
	We next recall that if $\mu$ is a $\sigma$-finite Borel measure on $M$ and $\mathbb{P}$ is a determinantal point process induced by a locally trace class operator $\Pi$ of orthogonal projection acting in the space $L_2(M,\mu)$, then the variance of an additive functional $S_f$, corresponding to a bounded measurable function $f$ of bounded support, is given by the formula 
\begin{equation}\label{var-f}
	\mathrm{Var} S_f = \displaystyle\frac{1}{2} \displaystyle\int\limits_{M}\displaystyle\int\limits_{M} |f(x)-f(y)|^2 \cdot |\Pi(x,y)|^2 d\mu(x) d\mu(y). 
\end{equation}
It therefore suffices, in order to establish the rigidity of the point process $\mathbb{P}_{\Pi}$, to find an increasing sequence of bounded Borel subsets $B^{(n)}$ exhausting $M$ and a sequence $f^{(n)}$ of bounded Borel functions of bounded support such that
$ f^{(n)}|_{B^{(n)}} \equiv 1$ and 
$$
	\lim\limits_{n\to\infty} \displaystyle\int\limits_{M}\displaystyle\int\limits_{M} |f^{(n)}(x)-f^{(n)}(y)|^2 |\Pi(x,y)|^2 d\mu(x) d\mu(y) = 0.
$$
\section{Rigidity in the continuous case.}
\subsection{Proof of Proposition \ref{gen-rig}.}

 Take $R>0$, $T>R$ and set 
$$
\varphi^{(R, T)}(x)=\begin{cases}
1-\displaystyle\frac{\log^+(|x|-R)}{\log(T-R)}  \  \mathrm{if} \  |x|\leq T;\\
0, |x|\geq T.\end{cases}
$$
To establish Proposition \ref{gen-rig}, it suffices to prove 
\begin{lemma}
If $\Pi$ satisfies the assumptions of Proposition \ref{gen-rig}, then, for any sufficiently large $R>0$, as $T\to\infty$ we have 
$
\Var_{\Prob_{\Pi}} S_{\varphi^{(R,T)}}\to 0.
$
\end{lemma}

Proof. We  estimate the double integral (\ref{var-f}) for the additive statistic $f=\varphi^{(R,T)}$. Of course, if $|x|, |y|<R$ or if $|x|, |y|>T$, then the expression under the integral sign is equal to zero. We will now estimate our integral 
over the domain  $$\{x,y\in \mathbb R^2: R<|x|,|y|<T\}$$ and complete the proof by estimating the smaller contribution of the 
domains $$\{x,y\in \mathbb R^2: 0<|x|<T<|y|\}, \ \{x,y\in \mathbb R^2: 0<|x|<R<|y|<\infty\}.$$
We consider these three cases separately. 

{\it The First Case}: $x,y\in \mathbb R^2: R<|x|,|y|<T$.

It is clear that for any $x,y$ satisfying $|x|, |y|>R$ there exists a constant $C(R)$ depending only on $R$ such that  we have  
$$
|\log^+(|x|-R)-\log^+(|y|-R)|\leq C(R)|\log |x|-\log |y||.
$$

Using the first assumption of Proposition\ref{gen-rig}, we now estimate the integral (\ref{var-f})  for the additive statistic $f=\varphi^{(R,T)}$  
 from above by the expression
\begin{equation}\label{sec-int}
\displaystyle\frac{\mathrm {const}}{(\log T)^{2}}\displaystyle\int\limits_R^T\displaystyle\int\limits_R^T \left(\displaystyle\frac{\log x -\log y}{x-y}\right)^2 \left(\displaystyle\frac{x^{2\alpha}}{y^{2\alpha}}+
\displaystyle\frac{y^{2\alpha}}{x^{2\alpha}}\right)
 dxdy,
\end{equation}
where the implied constant depends only on $R$. 
Introducing the variable $\la=y/x$ and recalling that $\alpha<1/2$,  we estimate the integral (\ref{sec-int}) from above by the expression
\begin{equation}\label{log-est}
\displaystyle\frac{\mathrm {const}}{(\log T)^{2}} \displaystyle\int\limits_R^T \displaystyle\frac {dx}{x} \displaystyle\int\limits_{T^{-1}}^T \left(\displaystyle\frac{\log \la}{\la-1}\right)^2(\la^{2\alpha}+\la^{-2\alpha})d\la=O\left(\log^{-1}T\right). 
\end{equation}

{\it The Second Case}: $x,y\in \mathbb R^2: |x|>R,|y|>T$. 

Next, we  consider the integral 
$$
\displaystyle\int\limits_R^T dx \displaystyle\int\limits_T^{\infty} (\varphi^{(R, T)}(x))^2 (\Pi(x,y))^2  dy,
$$
which (upon recalling that $x\leq y$ and making a scaling change of variable) can be estimated from above by the expression
$$
\displaystyle\frac{\mathrm{const}}{\log^2T}\displaystyle\int\limits_0^1 dx \displaystyle\int\limits_1^{\infty}\left(\frac{ y^{2\alpha}}{x^{2\alpha}}+1\right)\left(\displaystyle\frac{\log x}{x-y}\right)^2 
 dy=O(\log^{-2}T).
$$

{\it The Third Case.} $\{x,y\in \mathbb R^2: 0<|x|<R<|y|<\infty\}$.

Finally, we consider the integral 
$$
\displaystyle\int\limits_0^R dx  \displaystyle\int\limits_R^{\infty} (\varphi^{(R,T)}(y)-1)^2 (\Pi(x,y))^2  dy,
$$
in order to estimate which it suffices to estimate the integral 
$$
\displaystyle\int\limits_0^R dx  \displaystyle\int\limits_R^{\infty} (\log^+(y-R))^2 (\Pi(x,y))^2  dy
$$
which, using the second assumption of Proposition \ref{gen-rig}, we estimate from above by the expression
$$
\displaystyle\frac {{\mathrm{const}}}{\log^2T} \int_R^{\infty} \frac{(\log y)^2}{y^{1+2\varepsilon}} dy=O(\log^{-2}T). 
$$
where the implied constant, as always, depends only on $R$. 
The proposition is proved completely. 
\subsection{The case of integrable kernels.}
In applications, one often meets kernels  admitting an {\it integrable} representation
\begin{equation}\label{integ}
\Pi(x,y)=\displaystyle \frac{A(x)B(y) - B(x)A(y)}{x-y};
\end{equation}
with smooth functions $A$, $B$; 
the diagonal values of the kernel $\Pi$ are given by the formula
\begin{equation}\label{pi-diag}
\Pi(x, x) = A'(x)B(x)-A(x)B'(x).
\end{equation}

In this case, Proposition \ref{gen-rig} yields the following 
\begin{corollary}\label{rig-int}
If the kernel $\Pi$ admits an integrable representation (\ref{integ}) and, furthermore,  there exist $R>0$, $C>0$  and $\varepsilon>0$ such that 
\begin{enumerate}
\item for all $|x|<R$ we have $|A(x)|\leq C|x|^{-1/2+\varepsilon}$;  $|B(x)|\leq C|x|^{-1/2+\varepsilon}$;
\item for all $|x|>R$ we have $|A(x)|\leq C|x|^{1/2-\varepsilon}$;  $|B(x)|\leq C|x|^{1/2-\varepsilon}$,
\end{enumerate}
then the process $\Prob_{\Pi}$ is rigid.
\end{corollary}

Proof. Indeed, it is clear that both assumptions of Proposition \ref{gen-rig} are verified in this case. 

\section{Examples: the Bessel and the Airy Kernel.}
\subsection{The determinantal point process with the Bessel kernel.}
	
	 Take $s>-1$ and recall that the Bessel kernel is given by the formula 
$$
\mathbf {J}_s(x,y) = \displaystyle\frac{\sqrt{x}J_{s+1}(\sqrt{x})J_s(\sqrt{y})-\sqrt{y}J_{s+1}(\sqrt{y})J_s(\sqrt{x}) }{2(x-y)} , x,y > 0.
$$
By the Macchi-Soshnikov theorem, the Bessel kernel induces a determinantal point process $\Prob_{\mathbf{J}_s}$ on $\Conf({\mathbb R}_+)$.
\begin{proposition}\label{bes-rig}
The determinantal point process  $\Prob_{\mathbf{J}_s}$ is rigid. 
\end{proposition}

Proof. Indeed, this follows from Corollary \ref{rig-int}, the estimate $J_s(x)\sim x^{s/2}$, valid for small $x$ (cf. e.g. 9.1.10 in  in Abramowitz and Stegun \cite{AS})
 and the standard asymptotic expansion
$$
J_s(x)= \sqrt{\frac 2{\pi x}}\cos(x-s\pi/2-\pi/4)+O(x^{-1})
$$
of the Bessel function of a large argument (cf. e.g. 9.2.1 in Abramowitz and Stegun \cite{AS}).
Proposition \ref{bes-rig} is proved. 

\subsection{The determinantal point process with the Airy kernel.}

	Recall that the Airy kernel is given by the formula
$$
	\Ai(x,y) = \displaystyle\frac{\Ai(x)\Ai'(y) - \Ai(y)\Ai'(x)}{x-y},
$$
where
$$
	\Ai(x) = \displaystyle\frac{1}{\pi}\displaystyle\int\limits_0^{+\infty} \cos \left(\displaystyle\frac{t^3}{3}+xt\right)dt
$$
is the standard Airy function. 

By the Macchi-Soshnikov theorem, the Airy kernel infuces a determinantal point process $\Prob_{\Ai}$ on $\Conf({\mathbb R})$. 
In this case, we establish rigidity in the following slightly stronger form. 
\begin{proposition}\label{airy-rigid}
For any $D\in {\mathbb R}$, the random variable $\#_{(D, +\infty)}$ is measurable with respect to the $\Prob_{\Ai}$-completion of the sigma-algebra ${\mathcal F}_{(-\infty, D)}$.
\end{proposition}
 Proof. Again, we take  $R>0$, $T>R$ and set 
$$
\varphi^{(R,T)}(x)=\begin{cases} 
0, \mathrm{for} \  x<-T;\\
1-\displaystyle\frac{\log^+ (|x|-R)}{\log (T-R)},  \mathrm{for} \  -T<x<-R;\\
1, \ \mathrm{ for} \   x\geq -R. 
\end{cases}
$$
Since $\Prob_{\Ai}$-almost every trajectory admits only finitely many particles on the positive semi-axis, 
the additive functional $S_{\varphi^{(T)}}$ is  $\Prob_{\Ai}$-almost surely well-defined. It is immediate from (\ref{var-f}) that its variance is finite. 

\begin{lemma}\label{airy-var-est} For any fixed $R>0$, as $T\to\infty$, we have
$
\Var S_{\varphi^{(R,T)}}\to 0.
$
\end{lemma}
The proof of Lemma \ref{airy-var-est} is done in exactly the same way as that of Proposition \ref{gen-rig} and Corollary \ref{rig-int} , 
using standard power estimates for the Airy function and its derivative for negative values of the argument (cf. e.g. 10.4.60, 10.4.62 in Abramowitz and Stegun \cite{AS}) as well as the standard superexponential estimates  for the Airy function and its derivative for positive values of the argument (cf. e.g. 10.4.59, 10.4.61 in Abramowitz and Stegun \cite{AS}).
Proposition \ref{airy-rigid} follows immediately. 

\section{Rigidity of determinantal point processes with discrete phase space.}
\subsection{A general sufficient condition.}
Proposition \ref{gen-rig} admits a direct analogue in the case of a discrete phase space. 
\begin{proposition}\label{discr-rig}
Let $\Pi(x,y)$ be a kernel yielding an operator of orthogonal projection 
acting in  $L_2({\mathbb Z})$.
Assume that there exists $\alpha\in (0, 1/2)$, $\varepsilon>0$, and, for any $R>0$, a constant $C(R)>0$ such that   the following holds: 
\begin{enumerate}
\item if $|x|, |y|\geq R$, then 
$$
|\Pi(x,y)|\leq C(R)\cdot \frac{({x}/{y})^{\alpha}+ ({y}/{x})^{\alpha}}{|x-y|};
$$
\item if $|x|\leq R$, then  for all $y$ we have 
$$
\sum\limits_{x:|x|\leq R}|\Pi(x,y)|^2\leq \frac{C(R)}{1+y^{1+\varepsilon}}.
$$

\end{enumerate} 
Then the point process $\Prob_{\Pi}$ is rigid. 
\end{proposition}

The proof is exactly the same as that of Proposition \ref{gen-rig}.
The Corollary for integrable kernels assumes an even simpler form in the discrete case. 
\begin{corollary}\label{drig-int}
If the kernel $\Pi$ admits an integrable representation (\ref{integ}) and, furthermore,  there exist $R>0$, $C>0$  and $\varepsilon>0$ such that 
 for all $|x|>R$ we have $$|A(x)|\leq C|x|^{1/2-\varepsilon}, \  |B(x)|\leq C|x|^{1/2-\varepsilon},$$
then the process $\Prob_{\Pi}$ is rigid.
\end{corollary}

\subsection{The determinantal point process with the Gamma-kernel}
Let $\zp=1/2+{\mathbb Z}$ be the set of {\it half-integers}.
The Gamma-kernel with parameters $z, z^{\prime}$ is defined on $\zp\times \zp$ by the formula
\begin{multline}
\mathbf {\Gamma}_{z, \zr}(x,y)=\frac{\sin(\pi z)\sin(\pi\zr)}{\pi\sin(\pi(z-\zr))}\times \\
(\Gamma(x+z+1/2)\Gamma(x+\zr+1/2)\Gamma(y+z+1/2)\Gamma(y+\zr+1/2))^{-1/2}\times\\
\frac{\Gamma(x+z+1/2)\Gamma(y+\zr+1/2)-
\Gamma(x+\zr+1/2)\Gamma(y+z+1/2)}{x-y}.
\end{multline}
Following Borodin and Olshanski  (cf. Proposition 1.8 in \cite{BO-gamma}), we consider two cases: first, the case of the {\it principal} series, where $z^{\prime}={\overline z}\notin {\mathbb R}$ and 
the case of the complementary series, 
in which $z, \zr$ are real and, moreover, there exists an integer $m$ such that  $z, z^{\prime}\in (m, m+1)$.
In both these cases, the Gamma-kernel induces an operator of orthogonal projection acting in $L_2(\zp)$. We now establish the rigidity of the corresponding determinantal measure $\Prob_{\mathbf {\Gamma}_{z, \zr}}$ on $\Conf(\zp)$. 
We use Corollary \ref{drig-int}. In the case of the principal series, the functions $A, B$ giving the integrable representation for the Gamma-kernel, are  bounded above, so there is nothing to prove. In the case of the complementary series,  the standard asymptotics 
$$
\frac{\Gamma(x+z)}{\Gamma(x+\zr)}\sim x^{z-\zr}
$$
 (cf. e.g. 6.1.47 in  Abramowitz and Stegun \cite{AS})
allows us directly to apply Corollary \ref{rig-int} 
and thus to complete the proof of
\begin{proposition}
The determinantal point process with the Gamma-kernel is rigid for all values of the parameters $z$, $\zr$ belonging to the principal and the complementary series. 
\end{proposition}

{\bf Acknowledgements.} I am deeply grateful to Grigori Olshanski and Yanqi Qiu for useful discussions.
This work is supported by A*MIDEX project (No. ANR-11-IDEX-0001-02), financed by Programme ``Investissements d'Avenir'' of the Government of the French Republic managed by the French National Research Agency (ANR). It is also supported in part by the Grant MD-2859.2014.1 of the President of the Russian Federation, by the Programme ``Dynamical systems and mathematical control theory'' of the Presidium of the Russian Academy of Sciences, by the ANR under the project ``VALET'' (ANR-13-JS01-0010) of the Programme JCJC SIMI 1, and by the RFBR grant 13-01-12449.

\end{document}